\documentstyle[11pt]{article}
\newcommand{\R}{I \!  \! R}
\newcommand{\be}{\begin{equation}}
\newcommand{\ee}{\end{equation}}
\newcommand{\bea}{\begin{eqnarray}}
\newcommand{\eea}{\end{eqnarray}}

\newcommand{\la}{\label}
\newcommand{\xa}{\alpha}

\newcommand{\xd}{\delta}

\newcommand{\xe}{\varepsilon}
\newcommand{\e}{\epsilon}

\newcommand{\xl}{\lambda}

\newcommand{\xO}{\Omega}
\newcommand{\C}{{\bf (C)}}

\newcommand{\RR}{ (R)}

\newcommand{\ra}{\rightarrow}

\newtheorem{theorem}{Theorem}[section]
\newtheorem{lemma}[theorem]{Lemma}

\begin{document}

\title{\Large  \bf{ Critical  Hardy--Sobolev Inequalities} }
\newcommand{\affil}[1]{{\footnotesize\sl #1}}

\author{
 {\large S. Filippas,} \\
\affil{Department of Applied Mathematics,
         University of Crete,
         71409 Heraklion,  Greece} \\
\affil{and Institute of Applied and Computational Mathematics,
        FORTH, 71110 Heraklion, Greece} \\ \\
{\large V.  Maz'ya,} \\
\affil{Department of Mathematical Sciences, The University of
Liverpool, }\\
\affil{M\&O Building, Peach Street, Liverpool, L69 72L UK} \\
\affil{and Department of Mathematics, Ohio State University,
Columbus,OH 43210 USA} \\ \\
{\large \&  A. Tertikas} \\
\affil{Department of Mathematics,
         University of Crete,
         71409 Heraklion,  Greece} \\
\affil{and Institute of Applied and Computational Mathematics,
        FORTH, 71110 Heraklion, Greece}
}
\date{}

\maketitle

\begin{abstract}
We consider Hardy inequalities in $\R^n$, $n \geq 3$,
 with best constant that involve either
 distance to the boundary or  distance to a surface of
co-dimension $k<n$,  and we show that they can still be
 improved by adding a multiple of a whole range of critical norms
 that at the extreme case become precisely
the critical Sobolev norm.

\vskip 0.5\baselineskip

\centerline{\bf R\'esum\'e} \vskip 0.5\baselineskip \noindent Nous
\' etudions des inegalit\'es de Hardy dans $\R^n, n \geq 3$, avec
meilleure constante, li\'{e}e soit \`{a} la distance au bord, soit
\`{a} la distance \`a une surface de codimension $k<n$. Nous
obtenons des versions am\'eillor\'ees en ajoutant un certain
nombre  des normes critiques qui, au cas extr\`eme, sont
pr\'ecisement les normes de Sobolev critiques.

\smallskip

\end{abstract}

\section{Introduction}

 Let $\xO \subset  \R^n$ be a  domain and $K$ be a compact, $C^2$ manifold without boundary
embedded in $\R^n$, of co-dimension $k$,  $1 \leq k <n$. When
$k=1$ we assume that $K = \partial \xO$, whereas for $1< k <n$
 we assume that
 $K \cap \bar{\xO} \neq \emptyset$.
We  set $d(x) = {\rm dist}(x,K)$.

We also recall for $1<p$ and $p\neq k$ the following condition
that was introduced in
\cite{BFT},\\

\qquad~\qquad~~~~     ~~~~~~~~~~~~~~~$-\Delta_p
d^{\frac{p-k}{p-1}} \geq 0$~~ on~~
 $\Omega \setminus K$,  \hfill    \C
\\
where $\Delta_p$ is the p-Laplacian, that is  $\Delta_p u = {\rm
div} (|\nabla u|^{p-2}\nabla u).$ We note that for $k=1$,
condition $\C$ becomes
 $- \Delta d \geq 0,$
which is equivalent to the convexity of the domain  $\xO$ for
$n=2$, but it is a much weaker condition than  convexity  of $\xO$
for $n \geq 3$.

Under assumption $\C$ the following Hardy inequality holds true
\cite{BFT},
\be\la{1.2} \int_{\xO} |\nabla u|^p dx- \left| \frac{p-k}{p}
\right|^p
 \int_{\xO} \frac{|u|^p}{d^p} dx \geq 0, \qquad
 u \in C_0^{\infty}(\xO \setminus K),
\ee
where $\left| \frac{p-k}{p}
\right|^p$ is the best constant.

Here is our main result, which shows that inequality (\ref{1.2})
can  be improved by adding a multiple  of a whole range of
critical norms
 that at the extreme case become precisely
the critical Sobolev norm.
\begin{theorem}\la{nth1}
Let $2 \leq p <n$, $p \neq k <n$    and $p < q \leq
\frac{np}{n-p}$. Suppose that   $\xO \subset \R^n$ is a  bounded
domain  and  $K$ is  a compact, $C^2$ manifold without boundary
embedded in $\R^n$, of co-dimension $k$,  $1 \leq k <n$. When
$k=1$ we assume that $K=\partial \xO$, whereas for $1< k <n$
 we assume that
 $K \cap \bar{\xO} \neq \emptyset$.
\\
(i) If in addition $\xO$ and $K$  satisfy condition $\C$, then
 there exists a positive constant $c=c(\xO,K)$    such that
for all $ u \in C_0^{\infty}(\xO \setminus K)$, there holds
\be\la{1.3} \int_{\xO} |\nabla u|^p dx-  \left| \frac{p-k}{p}
\right|^p  \int_{\xO} \frac{|u|^p}{d^p} dx
  \geq
 c \left( \int_{\xO}  d^{-q+\frac{q-p}{p}n} |u|^{q} dx
 \right)^{\frac{p}{q}}.
\ee
 (ii) Without assuming condition $\C$, there exist a positive
 constant $c=c(n,k,p,q)$ independent of $\xO $, $K$ and a
 constant $M=M(\xO,K)$, such that
for all $ u \in C_0^{\infty}(\xO \setminus K)$, there holds
\be\la{1.4} \int_{\xO} |\nabla u|^p dx-  \left| \frac{p-k}{p}
\right|^p  \int_{\xO} \frac{|u|^p}{d^p} dx + M   \int_{\xO} |u|^p
dx
  \geq
 c \left( \int_{\xO}  d^{-q+\frac{q-p}{p}n} |u|^{q} dx
 \right)^{\frac{p}{q}}.
\ee
\end{theorem}

We note that the term in the right hand side of (\ref{1.3}) and
(\ref{1.4}) is optimal and in fact (\ref{1.3}) is a scale
invariant inequality. In the extreme case where
$q=\frac{np}{n-p}$, the term  in the  right hand side is precisely
the critical Sobolev term.

 The only  result  previously  known, in the spirit of estimate  (\ref{1.3}),
 concerns
 the  particular case where  $\xO =\R^n$,  $p=2$ and $K$ is affine, that is,
  $K=\{x \in \R^n~ |~x_1=x_2=\ldots=x_k=0 \}$,
 $1 \leq  k < n$, $k \neq 2$ and  has  been established in
\cite{M}. The case $p\neq 2$ was posed as an open question in
\cite{M}.
% (with $M=0$).

On the other hand the nonnegativity of the left hand side of
(\ref{1.4}) for $p=2$ has been shown in \cite{BM} for $K= \partial
\xO$. Other improvements  of the plain Hardy inequality involving
any arbitrary  subcritical $L^{q}$ term are presented in
\cite{FMT2} for the case where $\xO$ is a convex domain and $K=
\partial \xO$. For earlier results involving  improvements with some subcritical
$L^{q}$ terms see \cite{DD}.

We emphasize that in our Theorem the case $k=n$, which corresponds
to taking distance from  an interior  point, is excluded. As a
matter of fact estimate (\ref{1.3}) fails in this case.
 Indeed in this case, the optimal  improvement  of the plain Hardy inequality
  involves the critical Sobolev exponent, but
contrary to (\ref{1.3}) it also has a logarithmic correction \cite{FT}.

To establish Theorem \ref{nth1} a crucial step is to obtain
local estimates in  a neighborhood of  $K$, see
Theorem  \ref{thm5.1}.

 For other directions in improving Hardy inequalities we refer to
 \cite{AE}, \cite{BV}, \cite{BFT}, \cite{BM}, \cite{CM}, \cite{DELV},
\cite{GP}, \cite{HHL}, \cite{M}, \cite{MMP}, \cite{Ti1},
\cite{Ti2}, \cite{VZ} and references therein.

The paper is organized  as follows. In section 2  we  establish
auxiliary weighted Sobolev type inequalities, in the special case
where distance is taken from the boundary. We then use these
inequalities in section 3 to  derive Hardy--Sobolev inequalities
when distance is taken from the boundary. In sections 4 and 5 we
consider more general distance functions, where distance is taken
from a surface $K$  of co-dimension $k$, as well as other critical
norms via interpolation.

Some preliminary  results have  been  announced in \cite{FMT1}.

{\bf Acknowledgments}
 This work started when  VM
 visited    the University
of Crete, the  hospitality and support of which is   acknowledged.
 SF and AT acknowledge partial  support by the  RTN  European network
 Fronts--Singularities, HPRN-CT-2002-00274.
 AT acknowledges
the hospitality and support of Rutgers University.

\setcounter{equation}{0}

\section{Weighted inequalities  involving the distance function}
\la{vineq}

 Let     $\xO \subset \R^n$  be   a  bounded  domain with $C^2$ boundary
 and  $d(x) = {\rm dist}(x, \partial \xO)$.
 We denote by   $\xO_{\delta}:=\{x \in \xO: ~ {\rm dist}(x, \partial \xO) \leq \delta \}$
 a tubular neighborhood of $\partial \xO$, for $\delta$ small.
Then, for $\xd$ small we have that $d(x) \in C^2(\xO_{\delta})$.
Also, if  $x \in \xO_{\xd}$ approaches $x_0 \in \partial \xO \in C^2$
 then clearly $d(x) \ra  0$ and also
\[
 \Delta d(x) = (N-1) H(x_0) + O(d(x)),
\] where $H(x_0)$ is the mean curvature of $\partial \xO$
at $x_0$; see e.g., [GT] section 14.6.  As a consequence of this we have that
there exists a  $\xd^*$ sufficiently small and a positive constant
$c_0$  such that
\\

\qquad~\qquad~~~~$|d \Delta d | \leq c_0 d $,~~ in ~~~$\xO_{\xd}$,
~~~for  $0< \xd  \leq \xd^*$.
 \hfill    \RR
\\

\noindent
We say that a domain   $\xO \subset \R^n$ satisfies condition ($R$) if
there exists a $c_0$ and a  $\xd^*$ such that ($R$) holds. In case $d(x)$ is not
a $C^2$ function we interpret the inequality in  ($R$) in the weak sense,
that is
\[
|\int_{\xO_{\xd}}\, d \,  \Delta d \,  \phi \,  dx |
 \leq  c_0 \int_{\xO_{\xd}} \, d \, \phi \,  dx,
\qquad  \forall   \phi \in C^{\infty}_0(\xO), ~~\phi \geq 0.
\]

In our proofs, instead of assuming that $\xO$ is a bounded  domain of class
$C^2$ we will sometimes  assume that $\xO$ satisfies condition ($R$). Thus, some of
our results hold true for a larger class of domains.  For instance,
if  $\xO$ is a strip or an infinite cylinder, condition ($R$) is easily
seen to be satisfied even though $\xO$ is not bounded.

We first prove an $L^1$ estimate.
\begin{lemma}\la{lem.1}
Let $\xO$ be a  bounded domain  which satisfies condition ($R$). For any   $a>0$ and
$S  \in \left(0, \frac12  n \pi^{\frac12}[\Gamma(1+n/2)]^{-\frac1n} \right)$,
there exists  $\xd_0=\xd_0(a/c_0)$ such that for all $\xd \in (0, \xd_0]$ there holds
\be\la{12.3}
  \int_{\xO_{\xd}}
d^a |\nabla v| dx + \int_{\partial \xO_{\xd}^c} d^a  |v| dS_x ~ \geq ~
 S  \,  \| d^a v \|_{L^{\frac{N}{N-1}}(\xO_{ \xd})},
  ~~~~~~\forall v \in
C^{\infty}(\xO).
\ee
\end{lemma}

{\em Proof:}
We will  use  the following inequality:  If $V \subset  \R^n$ is any bounded  domain and
$u \in C^{\infty}(V)$ then
\be\la{12.5}
S_n \, \|u \|_{L^{\frac{n}{n-1}}(V)} \leq \|\nabla u \|_{L^1(V)} +
\|u\|_{L^1(\partial V)},
\ee
where  $S_n = n \pi^{\frac12}[\Gamma(1+n/2)]^{-\frac1n}$; see \cite{M}, p. 189.

For $V=\xO_{\xd}$ we apply (\ref{12.5}) to   $u=d^a v$, $v \in  C^{\infty}(\xO)$ to get
 \be\la{12.7}
S_n  \, \| d^a \, v \|_{L^{\frac{N}{N-1}}(\xO_{\xd})} \leq
 \int_{\xO_{\xd}}
d^a |\nabla v| dx + a  \int_{\xO_{\xd}}  d^{a-1}  |v| dx +
\int_{\partial \xO_{\xd}^c} d^a  |v| dS_x,
\ee
To estimate the middle term of the right hand side,
 noting that $\nabla d \cdot \nabla d=1$
a.e. and integrating by parts we have
\[
 a  \int_{\xO_{\xd}}  d^{a-1}  |v| dx =    \int_{\xO_{\xd}} \nabla  d^{a}  \cdot
 \nabla d  \,  |v| dx =  -  \int_{\xO_{\xd}} d^{a} \Delta d  |v| dx -
\int_{\xO_{\xd}} d^{a}  \nabla d \cdot \nabla  |v| dx +
\int_{\partial \xO_{\xd}^c}  d^{a} |v| dS_x
\]
Under our condition  ($R$) for $\xd$ small we have
$|d \Delta d| <c_0 d$ in $\xO_{\xd}$.  It follows  that
\be\la{12.9}
 (a-c_0 \xd)  \int_{\xO_{\xd}}  d^{a-1}  |v| dx  \leq
 \int_{\xO_{\xd}} d^{a} | \nabla v| dx +  \int_{\partial \xO_{\xd}^c}  d^{a} |v| dS_x.
\ee From (\ref{12.7}) and (\ref{12.9}) we get
\[
\frac{a-c_0 \xd}{2a - c_0 \xd}
S_n  \, \| d^a \, v \|_{L^{\frac{n}{n-1}}(\xO_{\xd})} \leq
 \int_{\xO_{\xd}}
d^a |\nabla v| dx +
\int_{\partial \xO_{\xd}^c} d^a  |v| dS_x.
\]
The result then
 follows  by taking
\be\la{12.10}
\xd_0 = \frac{a(S_n-2S)}{c_0(S_n-S)}.
\ee

$\hfill \Box$

We similarly have
\begin{lemma}\la{lem.2}
Let $\xO$ be a  domain  which satisfies condition ($R$).
 For any
$S  \in \left(0, \frac12 n v_n^{\frac1n} \right)$ and $a>0$
there exists  $\xd_0=\xd_0(a/c_0)$ such that for all $\xd \in (0, \xd_0]$ there holds
\be\la{12.11}
   \int_{\xO_{\xd}}
d^a |\nabla v| dx~
\geq
 ~S \| d^a v \|_{L^{\frac{n}{n-1}}(\xO_{ \xd})},
  ~~~~~~\forall v \in
C_0^{\infty}(\xO_{\xd}).
\ee
\end{lemma}
The proof is quite similar to that of the previous Lemma. Instead of (\ref{12.5})
one uses
 the
 $(p=1)$--Gagliardo--Nirenberg inequality valid for any $V \subset \R^n$,
and any  $u \in  C_0^{\infty}(V)$
\be\la{12.13}
\tilde{S}_n \, \|u \|_{L^{\frac{n}{n-1}}(V)} \leq \|\nabla u \|_{L^1(V)},
\ee
where  $\tilde{S}_n = n v_n^{\frac1n}$, and $ v_n$ denotes the volume of the unit ball
in $\R^n$.

We next prove
\begin{theorem}\la{thm1}
Let  $\xO$ be a  bounded domain of class $C^2$  and $1<p<n$. Then there exists a  $\xd_0=\xd_0(\xO,p,n)$
 such that for all $\xd \in (0, \xd_0]$ there holds
\be\la{12.15}
  \int_{\xO_{\xd}} d^{p-1} | \nabla v|^p dx +
\int_{\partial \xO_{\xd}^c} |v|^p dS_x ~  \geq ~
C(n,p) \|d^{\frac{p-1}{p}} v \|^p_{L^{\frac{np}{n-p}}(\xO_{\xd})},
 ~~~~~~\forall v \in
C^{\infty}(\xO),
\ee
with a constant $C(n,p)$ depending only on $n$ and $p$.
\end{theorem}

{\em Proof:} We will denote by  $C(p)$, $C(n,p)$  etc.   positive constants,
not necessarily the same in each occurrence, which depend {\em only} on their arguments.
 As a first step we will prove the following estimate:
\be\la{12.17}
C(n,p) \|d^{\frac{p-1}{p}} v \|^p_{L^{\frac{np}{n-p}}(\xO_{\xd})}
\leq  \int_{\xO_{\xd}} d^{p-1} | \nabla v|^p dx +
\| d^{\frac{p-1}{p}} v\|^p_{L^{\frac{(n-1)p}{n-p}}(\partial \xO_{\xd}^c)}.
\ee

To this end we apply estimate  (\ref{12.3}) to  $w=|v|^s$, $s= \frac{(n-1)p}{n-p}$
 with $a = \frac{(n-1)(p-1)}{n-p}>0$. Then,
\[
S(n,p) \left(
\int_{\xO_{\xd}} d^{\frac{n(p-1)}{n-p}}
|v|^{\frac{np}{n-p}} dx \right)^{\frac{n-1}{n}} \leq  s
\int_{\xO_{\xd}} d^{\frac{(n-1)(p-1)}{n-p}} |v|^{\frac{n(p-1)}{n-p}} |\nabla v|  dx
+  \int_{\partial \xO_{\xd}^c}
d^{\frac{(n-1)(p-1)}{n-p}} |v|^{\frac{(n-1)p}{n-p}} dS_x.
\]
We next  estimate the middle term
\bea
\int_{\xO_{\xd}} d^{\frac{(n-1)(p-1)}{n-p}} |v|^{\frac{n(p-1)}{n-p}} |\nabla v|  dx
& \leq  &
\left(
\int_{\xO_{\xd}} d^{\frac{n(p-1)}{n-p}}
|v|^{\frac{np}{n-p}} dx \right)^{\frac{p-1}{p}}
\left(
\int_{\xO_{\xd}} d^{p-1}
|\nabla v|^p dx \right)^{\frac{1}{p}}  \nonumber \\
& \leq  &
\e
\left(
\int_{\xO_{\xd}} d^{\frac{n(p-1)}{n-p}}
|v|^{\frac{np}{n-p}} dx \right)^{\frac{n-1}{n}} +c_{\e}
\left(
\int_{\xO_{\xd}} d^{p-1}
|\nabla v|^p dx \right)^{\frac{n-1}{n-p}},  \nonumber
\eea
whence,
\[
(S(n,p)-\e s)
\left(
\int_{\xO_{\xd}} d^{\frac{n(p-1)}{n-p}}
|v|^{\frac{np}{n-p}} dx \right)^{\frac{n-1}{n}} \leq
s c_{\e}
\left(
\int_{\xO_{\xd}} d^{p-1}
|\nabla v|^p dx \right)^{\frac{n-1}{n-p}}
+
\int_{\partial \xO_{\xd}^c}
d^{\frac{(n-1)(p-1)}{n-p}} |v|^{\frac{(n-1)p}{n-p}} dS_x.
\]
Raising the above estimate to the power  $\frac{n-p}{n-1}$ we easily
 obtain (\ref{12.17}).

To prove (\ref{12.15}) we need to combine (\ref{12.17}) with the following
estimate
\be\la{12.19}
C(n,p)
\| d^{\frac{p-1}{p}} v\|^p_{L^{\frac{(n-1)p}{n-p}}(\partial \xO_{\xd})}
%\left( \int_{\partial \xO_{\xd}}
%d^{\frac{(n-1)(p-1)}{n-p}} |v|^{\frac{(n-1)p}{n-p}} dx \right)^{\frac{n-p}{n-1}}
\leq
\int_{\xO_{\xd}} d^{p-1}
|\nabla v|^p dx +
\int_{\partial \xO_{\xd}^c} |v|^p dS_x.
\ee
In the rest of the proof we   will   show (\ref{12.19})
We note that the  norm  in the left hand side is the critical trace norm
of the function $ d^{\frac{p-1}{p}} v$. To estimate it we will  use the
critical trace inequality (\cite{b}, Proposition 1),
\be\la{12.21}
\| u\|^p_{L^{\frac{(n-1)p}{n-p}}(\partial \xO_{\xd})}
\leq C(n,p) \| \nabla u \|^p_{L^p(\xO_{\xd})} + M \|u \|^p_{L^p(\xO_{\xd})},
\ee
where $M=M(n,p, \xO)$ in general  depends on the domain $\xO$  as well.
For reasons that we will  explain later we will apply this estimate  not directly to
 $ d^{\frac{p-1}{p}} v$ but to the function $u=d^{\frac{p-1}{p} + \theta} v$ with
$\theta >0$ instead. More specifically we have
\bea
\| d^{\frac{p-1}{p}} v\|^p_{L^{\frac{(n-1)p}{n-p}}(\partial \xO_{\xd})} &  = &
\xd^{-\theta p} \, \,
\| d^{\frac{p-1}{p}+ \theta} v\|^p_{L^{\frac{(n-1)p}{n-p}}(\partial \xO_{\xd})}
\nonumber \\
& \leq &
\xd^{-\theta p} \, \,
\left( C(n,p) \| \nabla ( d^{\frac{p-1}{p}+ \theta} v)\|^p_{L^p(\xO_{\xd})}
+  M \|d^{\frac{p-1}{p} + \theta} v \|^p_{L^p(\xO_{\xd})} \right).  \nonumber
\eea
Now,
\[
 \| \nabla ( d^{\frac{p-1}{p}+ \theta} v)\|_{L^p(\xO_{\xd})} \leq
\left(\frac{p-1}{p} + \theta  \right) \|d^{-\frac{1}{p}+ \theta} v\|_{L^p(\xO_{\xd})} +
\|d^{\frac{p-1}{p} + \theta} \nabla v \|_{L^p(\xO_{\xd})},
\]
and
\[
\|d^{\frac{p-1}{p} + \theta} v \|_{L^p(\xO_{\xd})} \leq \xd \,
 \|d^{-\frac{1}{p}+ \theta} v\|_{L^p(\xO_{\xd})}.
\]
From the above three estimates we  conclude that \bea \|
d^{\frac{p-1}{p}} v\|^p_{L^{\frac{(n-1)p}{n-p}}(\partial
\xO_{\xd})} &  \leq &
C(p)  \xd^{-\theta p} \int_{\xO_{\xd}} d^{p-1 + p \theta} |\nabla v|^p dx  \nonumber \\
&  &  +
[C(n,p,\theta)+M \xd^{p}] \, \,
  \xd^{-\theta p}  \int_{\xO_{\xd}} d^{-1 + p \theta} |v|^p dx, \nonumber
\eea
whence, by  choosing $\xd$ sufficiently small,
\bea\la{12.23}
\| d^{\frac{p-1}{p}} v\|^p_{L^{\frac{(n-1)p}{n-p}}(\partial \xO_{\xd})} &  \leq &
C(p)  \xd^{-\theta p} \int_{\xO_{\xd}} d^{p-1 + p \theta} |\nabla v|^p dx  \nonumber \\
&  &  +
C(n,p,\theta) \, \,
  \xd^{-\theta p}  \int_{\xO_{\xd}} d^{-1 + p \theta} |v|^p dx.
\eea
To continue we will estimate the last term of the right hand side of (\ref{12.23}).
Consider the identity:
\be\la{12.23a}
  \theta p  d^{-1 + \theta p} =  - d^{\theta p} \Delta d + {\rm div}
(d^{\theta p} \nabla d)
\ee
We multiply it by $|v|^p$ and  integrate by parts  over
$\xO_{\xd}$ to get
\[
\theta p  \int_{\xO_{\xd}}  d^{-1 + \theta p} |v|^p dx = -  \int_{\xO_{\xd}}
d^{\theta p} \Delta d |v|^p dx  -p \int_{\xO_{\xd}} d^{\theta p} |v|^{p-1} \nabla d
 \cdot  \nabla |v|
 \, dx + \int_{\partial \xO_{\xd}^c} d^{\theta p} |v|^{p} dS_x.
\]
By our assumption ($R$) we have that $|d^{\theta p} \Delta d | \leq c_0 \xd \,  d^{-1 + \theta p}$.
On the other hand
\bea
 |p \int_{\xO_{\xd}} d^{\theta p} |v|^{p-1} \nabla d  \cdot  \nabla |v|
 dx | &   \leq  &  p \int_{\xO_{\xd}} d^{\theta p} |v|^{p-1} |\nabla v|  dx
\nonumber \\
 &   \leq  & p \e   \int_{\xO_{\xd}}  d^{-1 + \theta p} |v|^p dx
+ p  c_{\e}  \int_{\xO_{\xd}} d^{p-1 + p \theta} |\nabla v|^p dx.
\nonumber
\eea
Putting together the last estimates we get
\be\la{12.24}
(\theta p-c_0 \xd -p \e)  \int_{\xO_{\xd}}  d^{-1 + \theta p} |v|^p dx \leq
p c_{\e}  \int_{\xO_{\xd}} d^{p-1 + p \theta} |\nabla v|^p dx
 + \int_{\partial \xO_{\xd}^c} d^{\theta p} |v|^{p} dS_x,
\ee
whence, choosing $\xd$, $\e$ sufficiently small,
\be\la{12.25}
C(p,\theta) \int_{\xO_{\xd}}  d^{-1 +p \theta } |v|^p dx \leq
C(p)  \int_{\xO_{\xd}} d^{p-1 + p \theta} |\nabla v|^p dx
 + \int_{\partial \xO_{\xd}^c} d^{p \theta } |v|^{p} dS_x.
\ee
Combining   (\ref{12.23})  and (\ref{12.25}) we obtain
\bea\la{12.27}
C(n, p, \theta) \| d^{\frac{p-1}{p}} v\|^p_{L^{\frac{(n-1)p}{n-p}}(\partial \xO_{\xd})}
 & \leq  &
 \xd^{-\theta p}  \int_{\xO_{\xd}} d^{p-1 + p \theta} |\nabla v|^p dx +
 \xd^{-\theta p} \int_{\partial \xO_{\xd}^c} d^{p \theta } |v|^{p} dS_x \nonumber \\
 & \leq  &
  \int_{\xO_{\xd}} d^{p-1} |\nabla v|^p dx +
 \int_{\partial \xO_{\xd}^c}  |v|^{p} dS_x.
\eea
By choosing a  specific  value of $\theta$, e.g., $\theta=1$,   we get (\ref{12.19}).
We  note that estimate (\ref{12.25}) fails if $\theta=0$, and this is the reason
for introducing this artificial parameter.

$\hfill \Box$

We next  have
\begin{theorem}\la{thm2}
Let  $\xO  \subset   \R^n $ be a  domain satisfying ($R$)  and $1<p<n$.
 Then   there exists a  $\xd_0=\xd_0(c_0,p,n)$
 such that for all $\xd \in (0, \xd_0]$ there holds
\be\la{12.30}
 \int_{\xO_{\xd}} d^{p-1} | \nabla v|^p dx~
\geq ~
C(n,p) \|d^{\frac{p-1}{p}} v \|^p_{L^{\frac{np}{n-p}}(\xO_{\xd})},
 ~~~~~~\forall v \in
C_0^{\infty}(\xO_{\xd}),
\ee
with a constant $C(n,p)$ depending only on $n$ and $p$.
\end{theorem}

{\em Proof:} One works  as in the derivation of (\ref{12.17}),
using however  (\ref{12.11}) in the place of (\ref{12.3}). We omit the details.

We finally   establish the following:
\begin{theorem}\la{thm3}
Let $1 < p <n$ and $D = \sup_{x \in \xO} d(x)< \infty$. We assume
that  $\xO$ is  a   domain satisfying both conditions \C ~and
($R$). Then   there exists a positive  constant $C=C(n,p,c_0 D)$
such that for any
 $ v \in  C_0^{\infty}(\xO)$,
\be\label{12.40}
\int_{\xO} d^{p-1} |\nabla v|^p dx +
\int_{\xO} (-\Delta d) |v|^p   dx~
\geq~
 C  \| d^{\frac{p-1}{p}} v \|^p_{L^{\frac{np}{n-p}}(\xO)}.
\ee
%where $C$ depends on $n$, $p$, $\xe$ and  $\frac{\xd_0}{D}$.
\end{theorem}

{\em Proof:}  We first define suitable cutoff functions supported near the boundary.
Let $\xa(t) \in C^{\infty}([0,\infty))$ be a nondecreasing function such that
 $\xa(t) = 1$ for
$t \in [0,1/2)$,  $\xa(t) = 0$ for  $t \geq 1$ and $|\xa'(t)| \leq C_0$.
 For $\xd$ small we define  $\phi_{\xd}(x):=\xa(\frac{d(x)}{\xd}) \in C^2_0(\xO)$. Note that
 $\phi_{\xd}=1$ on $\xO_{\delta/2}$,
 $\phi_{\xd}=0$ on $\xO_{ \delta}^c$ and $| \nabla \phi_{\xd}|=|\xa' (\frac{d(x)}{\xd}) |
\frac{|\nabla d(x)|}{\xd}  \leq   \frac{C_0}{\xd}$ with  $C_0$ a universal constant.

For $v \in  C_0^{\infty}(\xO)$ we write  $v= \phi_{\xd} v + (1-\phi_{\xd}) v$.
The function $\phi_{\xd} v$ is compactly supported in $\xO_{\xd}$, and by Lemma \ref{lem.2} we
have
\be\la{12.41}
S  \, \| d^a \,\phi_{\xd} v \|_{L^{\frac{n}{n-1}}(\xO_{ \xd})} \leq
\int_{\xO} d^{a} |\nabla (\phi_{\xd} v)| dx.
\ee
On the other hand $(1-\phi_{\xd}) v$  is compactly supported in $\xO_{\xd/2}^c$ and using
(\ref{12.13}) we have
\be\la{12.43}
C(n) \, \| d^a \,(1-\phi_{\xd}) v \|_{L^{\frac{n}{n-1}}(\xO)} \leq
\left(\frac{2D}{\xd} \right)^a
\int_{\xO} d^{a} |\nabla ((1-\phi_{\xd}) v)| dx.
\ee
Combining (\ref{12.41}) and (\ref{12.43}) and using elementary estimates, we obtain
the following $L^1$ estimate:
\be\la{12.45}
 C(a,n, \frac{\xd}{D})\, \,  \| d^a v \|_{L^{\frac{n}{n-1}}(\xO)} \leq
 \int_{\xO}  | d^a  \nabla v | dx +
\int_{\xO_{ \delta} \setminus \xO_{ \delta/2}} d^{a-1}  |v| dx.
\ee
We next derive the corresponding $L^p$, $p>1$ estimate. To this end
 we replace $v$ by
 $|v|^s$ with  $s=\frac{p(n-1)}{n-p}$  in (\ref{12.45})  to obtain
\bea
   C(a,n,p, \frac{\xd}{D})\, \,  \left( \int_{\xO}
d^{\frac{an}{n-1}} |v|^{\frac{np}{n-p}} dx \right)^{\frac{n-1}{n}}
& \leq &
 s  \int_{\xO}  d^a |v|^{\frac{n(p-1)}{n-p} } |\nabla v| dx  \nonumber \\
&  & +   \int_{\xO_{ \delta} \setminus \xO_{ \delta/2}}
d^{a-1}|v|^{1+\frac{n(p-1)}{n-p} }dx. \nonumber
\eea
Using Holders inequality in both terms of the right hand side of this we get after simplifying,
\bea\la{12.47}
   C(a,n,p, \frac{\xd}{D})\,  \left( \int_{\xO} d^{\frac{an}{n-1}}
 |v|^{\frac{np}{n-p}} dx \right)^{\frac{n-p}{np}}
& \leq &  s   \left( \int_{\xO} d^{\frac{a(n-p)}{n-1}}|\nabla v|^p  \right)^{1/p}
\nonumber \\
& & +   \left( \int_{\xO_{ \delta} \setminus \xO_{ \delta/2}}
 d^{\frac{a(n-p)}{n-1}-p} |v|^p \right)^{1/p}.
\eea
For
 $a=\frac{(n-1)(p-1)}{n-p}>0$, this yields
\be\la{12.49}
  C(n,p, \frac{\xd}{D})\,   \| d^{\frac{p-1}{p}} v \|^p_{L^{\frac{np}{n-p}}(\xO)} \leq
\int_{\xO} d^{p-1} |\nabla v|^p dx +
 \int_{\xO_{ \delta} \setminus \xO_{ \delta/2}}  d^{-1} |v|^p dx.
\ee
 We note that condition $\C$
has not been used so far and therefore all previous estimates
 are valid even for general domains.

To complete the proof we will  estimate the last term in (\ref{12.49}).
For $\theta>0$, we clearly have
\be\la{12.50}
\left( \frac{\xd}{2} \right)^{p \theta}
 \int_{\xO_{ \delta} \setminus \xO_{ \delta/2}}  d^{-1} |v|^p dx
\leq
\int_{\xO_{ \delta} \setminus \xO_{ \delta/2}}  d^{-1+p \theta} |v|^p dx
\leq
\int_{\xO} d^{-1+p \theta} |v|^p dx.
\ee
To estimate the last term we work as in (\ref{12.23a})--(\ref{12.25}).
Thus, we start from the identity (\ref{12.23a}), multiply by $|v|^p$ and
 integrate
by parts in $\xO$. Now there are no boundary terms and also  the
term containing $\Delta d$ is  not a lower order term anymore and has
 to be kept.
 Notice however that  because of condition $\C$
 we have that $-\Delta d  \geq 0$ in the  distributional sense.
Without reproducing the details we write the analogue of
(\ref{12.25}) which is \be\la{12.51} C(p,\theta) \int_{\xO      }
d^{-1 +p \theta } |v|^p dx \leq C(p)  \int_{\xO}       d^{p-1 + p
\theta} |\nabla v|^p dx
 + \int_{\xO}   d^{p \theta }(-\Delta d) |v|^{p} dx.
\ee
Combining (\ref{12.50}) and (\ref{12.51}) and recalling that $d \leq D$, we get
\be\la{12.52}
 C(p, \theta) \left( \frac{\xd}{D} \right)^{ p \theta} \,
 \int_{\xO_{ \delta} \setminus \xO_{ \delta/2}}  d^{-1} |v|^p dx
\leq
 \int_{\xO}       d^{p-1} |\nabla v|^p dx
 + \int_{\xO}   (-\Delta d) |v|^{p} dx.
\ee
Choosing e.g., $\theta=1$ and combining (\ref{12.52}) and (\ref{12.49})
 the result
follows.
 The dependence of the constant $C$ in  (\ref{12.40}) on the
domain $\xO$ enters through the ratio $\xd/D$.  By Lemma \ref{lem.2}
(cf  (\ref{12.10}))  we obtain that the dependence of $C$ on $\xO$ enters
through $c_0 D$. We also note  that $C(n,p, \infty)=0$.

$\hfill \Box$

%%%%%%%%%%%%%%%%%%%%%%%%%%%%%%%%%%%%%%%%%%%

\setcounter{equation}{0}
\section{Hardy-- Sobolev inequalities}

%%%%%%%%%%%%%%%%%%%%%%%%%%%%%%%%%%%%%%%%%
Here we will prove various Hardy Sobolev inequalities.
 Let $d(x) = {\rm dist}(x, \partial \xO)$ and $V \subset \xO$.
For $p>1$,  and  $u \in C_{0}^{\infty}(\xO)$ we set
\be\la{3.defi}
I_{p}[u](V):=\int_{V}  |\nabla u|^pdx-
\left(\frac{p-1}{p}\right)^p  \int_{V}
\frac{|u|^p}{d^{p}}dx.
\ee
For simplicity we  also write $I_{p}[u]$ instead of $I_{p}[u](\xO)$.
We next put
\be\la{3.3}
u(x)=d^{\frac{p-1}{p}}(x) v(x).
\ee
We first prove an auxiliary inequality
\begin{lemma}\la{lem3.1}
 For $p \geq 2$, there exists positive constant $c=c(p)$ such that
\be\la{3.5}
I_{p}[u](V) \geq c(p) \int_{V} d^{p-1} |\nabla v|^p dx  +
\left(\frac{p-1}{p}\right)^{p-1}
 \int_{V} \nabla d \cdot \nabla |v|^p dx.
%\left(\frac{p-1}{p}\right)^{p-1}  \int_{\xO} (-\Delta d) |v|^p dx.
\ee
%In particular, for  $p=2$ we have the {\it identity:}
%\be\la{3.6}
%I_2[u] =  \int_{\xO} d |\nabla v|^2 dx  +
%\frac12  \int_{\xO} (-\Delta d) |v|^2 dx.
%\ee
\end{lemma}
{\em Proof:}
We have that
\[
\nabla u =\frac{p-1}{p} d^{\frac{p-1}{p}-1}
v \nabla d + d^{\frac{p-1}{p}} \nabla v =: a + b.
\]
For $p \geq 2 $ we have that for $a$, $b \in \R^n$,
\[
|a+b|^p - |a|^p \geq c(p)|b|^p + p |a|^{p-2} a \cdot b.
\]
Using this we obtain
\be\la{1.10}
I_{p}[u](V) \geq c(p) \int_{V}d^{p-1} |\nabla v|^p  dx +
 \left(\frac{p-1}{p}\right)^{p-1}
 \int_{V} \nabla d \cdot \nabla |v|^p dx.
\ee
which is the sought for estimate.

 $\hfill \Box$

We first establish estimates in $\xO_{\xd}$.
\begin{theorem}\la{thm3.1}
Let $2 \leq p <n$.  We   assume that  $\xO$ is a bounded domain of class $C^2$.
 Then, there exists a $\xd_0=\xd_0(p,n,\xO)$ such that  for $0< \xd \leq \xd_0 $
and all  $ u \in C_0^{\infty}(\xO)$
\be\la{3.10}
\int_{\xO_{\xd}} |\nabla u|^p dx- \left( \frac{p-1}{p}\right)^{p}
 \int_{\xO_{\xd}} \frac{|u|^p}{d^p} dx
\geq C
 \left( \int_{\xO_{\xd}} |u|^{\frac{np}{n-p}} dx \right)^{\frac{n-p}{n}},
\ee
where  $C=C(n,p)>0$ depends only on $n$ and $p$.
\end{theorem}

{\em Proof:} Using   Lemma \ref{lem3.1}  we have that
\[
C(p) \, I_{p}[u](\xO_{\xd}) \geq  \int_{\xO_{\xd}}  d^{p-1} |\nabla v|^p  dx +
 \int_{\xO_{\xd}} \nabla d \cdot \nabla |v|^p dx.
\]
Integrating by parts the last term we get
\be\la{3.11}
C(p) \, I_{p}[u](\xO_{\xd}) \geq  \int_{\xO_{\xd}}  d^{p-1} |\nabla v|^p  dx +
 \int_{\xO_{\xd}}(-\Delta d)  |v|^p dx +  \int_{\partial \xO_{\xd}^c}
 |v|^p dS_x.
\ee
We next estimate the middle term of the right hand side.
By  condition ($R$)  we have
\be\la{3.13}
| \int_{\xO_{\xd}}(-\Delta d)  |v|^p dx| \leq c_0 \int_{\xO_{\xd}}  |v|^p dx.
\ee
 Starting from the identity $1 + d \Delta d = {\rm div} (d \nabla d)$, we multiply it  by $|v|^p$ and
integrate by parts over $\xO_{\xd}$
to get
\[
\int_{\xO_{\xd}}|v|^p dx +
\int_{\xO_{\xd}} d \Delta d |v|^p dx  = -p \int_{\xO_{\xd}} d |v|^{p-1}
 \nabla d  \cdot \nabla |v| dx +  \xd
\int_{\partial \xO_{\xd}^c} |u|^p  dS.
\]
Using once more  ($R$)
and standard
inequalities we get
\[
(1- \xd c_0 - \xe p)
\int_{\xO_{\xd}}|v|^p dx  \leq
 \xd p C_{\xe}
\int_{\xO_{\xd}}d^{p-1} |\nabla v|^p dx
+   \xd \int_{\partial \xO_{\xd}^c} |u|^p  dS,
\]
whence for $\xe$, $\xd$ sufficiently small,
\be\la{3.15}
\int_{\xO_{ \xd}}|v|^p dx  \leq
C(p)  \xd \int_{\xO_{ \xd}}d^{p-1} |\nabla v|^p dx
+ C(p)  \xd \int_{\partial \xO_{\xd}^c} |u|^p  dS.
\ee
Combining (\ref{3.11}), (\ref{3.13}) and  (\ref{3.15}) we obtain,
\be\la{3.17}
C(p) \, I_{p}[u](\xO_{\xd}) \geq  \int_{\xO_{\xd}}  d^{p-1} |\nabla v|^p  dx +
 \int_{\partial \xO_{\xd}^c}
 |v|^p dS_x.
\ee
To complete the proof we now use Theorem \ref{thm1}, that is,
\bea\la{3.19}
 \int_{\xO_{\xd}} d^{p-1} | \nabla v|^p dx +
\int_{\partial \xO_{\xd}^c} |v|^p dS_x & \geq &
C(n,p) \|d^{\frac{p-1}{p}} v \|^p_{L^{\frac{np}{n-p}}(\xO_{\xd})}
\nonumber \\
& = &  C(n,p) \| u  \|^p_{L^{\frac{np}{n-p}}(\xO_{\xd})}.
\eea
The result then follows from  (\ref{3.17}) and  (\ref{3.19})

$\hfill \Box$

Next we   prove:
%\ref{thm3.5}. More precisely we have
\begin{theorem}\la{thm3.2}
Let $2 \leq p <n$.  We   assume that  $\xO$ is a bounded domain of class $C^2$.
Then
 there exists  positive constants $M=M(n,p,\xO)$ and  $C=C(n,p)$ such that
for all $ u \in C_0^{\infty}(\xO)$,
there holds
\be\la{3.21}
\int_{\xO} |\nabla u|^p dx-  \left( \frac{p-1}{p} \right)^p  \int_{\xO} \frac{|u|^p}{d^p} dx
 + M   \int_{\xO} |u|^p dx   \geq
 C \left( \int_{\xO} |u|^{\frac{np}{n-p}} dx \right)^{\frac{n-p}{n}}.
\ee
We emphasize that $C(n,p)$ is independent of $\xO$.
\end{theorem}

{\em Proof:} Clearly we have
\be\la{3.23}
I_{p}[u](\xO)=I_{p}[u](\xO_{\xd})+ I_{p}[u](\xO_{\xd}^c).
\ee
By Theorem    \ref{thm3.1}    for $\xd$ small we have
\be\la{3.25}
I_{p}[u](\xO_{\xd}) \geq C(n,p) \| u  \|^p_{L^{\frac{np}{n-p}}(\xO_{\xd})}.
\ee
Since $d(x) \geq \xd$ in $\xO_{\xd}^c$,
\be\la{3.27}
I_{p}[u](\xO_{\xd}^c) \geq \int_{\xO_{\xd}^c} |\nabla u|^p dx -
 \left( \frac{p-1}{p \xd} \right)^p
 \int_{\xO_{\xd}^c} |u|^p dx.
\ee
Using the Sobolev embedding of $L^{\frac{np}{n-p}}(\xO_{\xd}^c)$ into
 $W^{1,p}(\xO_{\xd}^c)$,  see \cite{H}, Theorem 4.1, we get
\[
 \| u \|^p_{L^{\frac{np}{n-p}}(\xO_{\xd}^c)}  \leq
C(n,p) \int_{\xO_{\xd}^c} |\nabla u|^p dx + C(n,p, \xO) \int_{\xO_{\xd}^c} |u|^p dx.
\]
From this and (\ref{3.27}) we get \be\la{3.29}
I_{p}[u](\xO_{\xd}^c) \geq C(n,p)  \| u
\|^p_{L^{\frac{np}{n-p}}(\xO_{\xd}^c)} - C(n,p, \xO) \int_{\xO}
|u|^p dx. \ee The result follows from (\ref{3.23}), (\ref{3.27})
and (\ref{3.29}).

$\hfill \Box$

%\subsection{Under convexity}

We finally   show
\begin{theorem}\la{thm3.3}
Let $2 \leq  p <n$ and $D = \sup_{x \in \xO} d(x)< \infty$. We
assume that  $\xO$ is  a   domain satisfying both conditions $\C$
and ($R$). Then   there exists a positive  constant $C=C(n,p,c_0
D)$ such that for any
  $ u \in C_0^{\infty}(\xO)$
there holds
\be\la{3.31}
\int_{\xO} |\nabla u|^p dx- \left( \frac{p-1}{p}\right)^{p} \int_{\xO} \frac{|u|^p}{d^p} dx
\geq C \left( \int_{\xO} |u|^{\frac{np}{n-p}} dx \right)^{\frac{n-p}{n}}.
\ee
\end{theorem}

{\em Proof:} Working as in the derivation of (\ref{3.11}) we get
\[
C(p) \, I_{p}[u](\xO) \geq  \int_{\xO}  d^{p-1} |\nabla v|^p  dx +
 \int_{\xO}(-\Delta d)  |v|^p dx.
\]
The result  then follows  from Theorem \ref{thm3}.

 $\hfill \Box$

\setcounter{equation}{0}
\section{Extensions}

Here we will extend the previous inequalities in two directions.
First by considering different distant functions and secondly by
 interpolating between the Sobolev
 $L^{\frac{pn}{n-p}}$   norm    and the $L^p$ norm. This way we will obtain
new scale invariant inequalities.

We  denote by  $K$ a surface
 embedded in
$\R^n$, of codimension $k$,  $1 < k <n$.
We also allow for the extreme cases  $k=n$ or 1, with the following
convention. In  case $k=n$,  $K$ is identified with the origin,
that is $K=\{0\}$, assumed to be in the interior of $\xO$.
In case $k=1$, $K$ is identified with $\partial \xO$.

From now on distance is taken from $K$, that is,
 $d(x) = {\rm dist}(x, K)$.
 We also set  $K_{\delta}:=\{x \in \xO: ~ {\rm dist}(x, K) \leq \delta \}$
 is  a tubular neighborhood of $K$, for $\delta$ small,
and  $K_{ \delta}^{c} :=\xO \setminus K_{\xd}$.

We say that $K$ satisfies condition ($R$) whenever
there exists a  $\xd^*$ sufficiently small and a positive constant
$c_0$  such that
\\

\qquad~\qquad~~~~$|d \Delta d+1-k | \leq c_0 d $,~~ in ~~~$K_{\xd}$,
~~~for  $0< \xd  \leq \xd^*$;
 \hfill    \RR
\\

\noindent
 For $k=1$ this
coincides with condition ($R$) of section 2. For  $k>1$,  if $K$ is
a compact, $C^2$ surface  without boundary, then condition  ($R$) is satisfied;
see,  e.g.,  \cite{AS} Theorem 3.2 or \cite{S} section 3.

%We assume that
%\be\la{14.2}
% D :=\sup_{x \in \xO}d(x)< \infty.
%\ee

We next present an interpolation Lemma.
\begin{lemma}
\la{interlem}
 Let  $a$, $b$, $p$ and $q$
be such that
\be\la{cond}
 1 \leq p<n, \quad    p < q \leq \frac{pn}{n-p}, \quad {\rm and}
  \quad b=a-1 + \frac{q-p}{qp}n.
\ee
 Then
for any $\eta>0$, there holds
\be\la{4.11}
 \|d^b v\|_{L^q(\xO)} \leq \xl  \eta^{-\frac{1-\lambda}{\lambda}}
 \|d^a v\|_{L^{\frac{pn}{n-p}}(\xO)} + (1-\xl) \eta
\|d^{a-1}v\|_{L^p(\xO)},  \qquad  \forall ~v \in C^{\infty}(\xO),
\ee
where
\be\la{4.13}
0 < \lambda := \frac{n(q-p)}{qp} \leq 1.
\ee
\end{lemma}
{\em Proof:} For $p_s:= \frac{pn}{n-p}$  and $\lambda$ as  in  (\ref{4.13})
 we   use  Holder's inequality to obtain
\bea
\int_{\xO} d^{qb} |v|^{q} dx &  = & \int_{\xO} (d^{a \lambda q} |v|^{\lambda q})
(d^{q(b-a \lambda)} |v|^{q(1-\lambda)}) dx  \nonumber \\
& \leq &  \left( \int_{\xO} d^{ap_s} |v|^{p_s} dx
\right)^{\frac{\lambda q }{p_s}}
 \left( \int_{\xO} d^{p(a-1)}|v|^p dx \right)^{\frac{(1-\lambda) q }{p}},
   \nonumber
\eea
that is,
\[
 \|d^b v\|_{L^q(\xO)} \leq  \|d^a v\|^{\lambda }_{L^{\frac{pn}{n-p}}(\xO)}~
\|d^{a-1} v\|^{1-\lambda}_{L^p(\xO)}~.
\]
Combining this with  Young's inequality
\be\la{4.14}
X^{\lambda } Y^{1-\lambda}  \leq  \xl \eta^{-\frac{1-\lambda}{\lambda}} X + (1-\xl)
\eta Y, \qquad
\eta>0,
\ee
the result follows.
$\hfill \Box$

We first  prove inequalities in $K_{\xd}$.
\begin{lemma}\la{lem.41}
Let $\xO \subset \R^n$ be a  bounded  domain
 and $K$ a $C^2$ surface of codimension $k$,
satisfying condition ($R$). We also assume that
\be\la{4.15}
   p=1 < q \leq \frac{n}{n-1},~~~ b=a-1+ \frac{q-1}{q}n, ~~~
{\rm and}~~~ a \neq 1-k.
\ee
Then there exists  a  $\xd_0=\xd_0(\frac{|a+k-1|}{c_0})$
 and  $ C=C(a,q,n,k)>0$
such that for all $\xd \in (0, \xd_0]$ there holds
\be\la{4.17}
  \int_{K_{\xd}}
d^a |\nabla v| dx + \int_{\partial K_{\xd}} d^a  |v| dS_x ~
\geq ~
 C \| d^b v \|_{L^{q}(K_{ \xd})},
~~~~~~\forall v \in
C_0^{\infty}(\xO \setminus K).
\ee
\end{lemma}
{\em Proof:}
 Using the interpolation inequality
 (\ref{4.11}) in  $K_{ \xd}$  with $\eta=1$  we get
\bea\la{4.19}
\|d^b v\|_{L^q(K_{ \xd})} & \leq  &
\frac{n(q-1)}{q}
 \|d^a v\|_{L^{\frac{N}{N-1}}(K_{ \xd})} +
\frac{q-n(q-1)}{q}
\|d^{a-1}v\|_{L^1(K_{ \xd})}  \nonumber \\
 & \leq  & C(n,q) \left(\|d^a v\|_{L^{\frac{N}{N-1}}(K_{ \xd})} +\int_{K_{\xd}}
d^{a-1} | v| dx \right).
\eea
For $V=K_{\xd}$ we apply (\ref{12.5}) to   $u=d^a v$, $v \in  C^{\infty}(\xO)$ to get
 \be\la{4.21}
S_n  \, \| d^a \, v \|_{L^{\frac{n}{n-1}}(K_{\xd})} \leq
 \int_{K_{\xd}}
d^a |\nabla v| dx + |a|  \int_{K_{\xd}}  d^{a-1}  |v| dx +
\int_{\partial K_{\xd}} d^a  |v| dS_x.
\ee
Combining (\ref{4.19}) and (\ref{4.21}) we get  the analogue of (\ref{12.7}) which is
 \be\la{4.22}
C(a,q,n)  \, \| d^b \, v \|_{L^{q}(K_{\xd})} \leq
 \int_{K_{\xd}}
d^a |\nabla v| dx +   \int_{K_{\xd}}  d^{a-1}  |v| dx +
\int_{\partial K_{\xd}} d^a  |v| dS_x.
\ee
It remains to  estimate the middle term of the right hand side. Noting that
 $\nabla d \cdot \nabla d=1$
a.e. and integrating by parts in  $K_{\xd}$  we have
\[
 a  \int_{K_{\xd}}  d^{a-1}  |v| dx =    \int_{K_{\xd}} \nabla  d^{a}  \cdot
 \nabla d  \,  |v| dx =  -  \int_{K_{\xd}} d^{a} \Delta d  |v| dx -
\int_{K_{\xd}} d^{a}  \nabla d \cdot \nabla  |v| dx +
\int_{\partial K_{\xd}}  d^{a} |v| dS_x,
\]
whence,
\[
(a+k-1) \int_{K_{\xd}}  d^{a-1}  |v| dx - \int_{K_{\xd}} d^{a-1} (d \Delta d +1-k)  |v| dx -
\int_{K_{\xd}} d^{a}  \nabla d \cdot \nabla  |v| dx +
\int_{\partial K_{\xd}}  d^{a} |v| dS_x.
\]
Using ($R$)  we easily arrive at the analogue of (\ref{12.9}),
that is,
\be\la{4.23}
(|a+k-1|-c_0 \xd)\int_{K_{\xd}}  d^{a-1}  |v| dx \leq
\int_{K_{\xd}} d^{a} | \nabla  v| dx +
\int_{\partial K_{\xd}}  d^{a} |v| dS_x.
\ee
For estimate (\ref{4.23}) to be useful we need $|a+k-1|>0$, whence the
restriction $a \neq 1-k$.
The result then follows from (\ref{4.22}) and (\ref{4.23}), taking
e.g.,  $\xd_0 = \frac{|a+k-1|}{2c_0}$.

$\hfill \Box$

We next present the analogue of Lemma \ref{lem.2}
\begin{lemma}\la{lem.42}
Let $\xO \subset \R^n$ be a domain  and $K$ a  surface of co-dimension $k$,
satisfying condition ($R$). We also assume
\[
   p=1 < q \leq \frac{n}{n-1},~~~ b=a-1+ \frac{q-1}{q}n, ~~~
{\rm and}~~~ a \neq 1-k.
\]
Then, there exists a   $\xd_0=\xd_0(\frac{|a+k-1|}{c_0})$
and a $C= C(a,q,n,k)>0$,  such that
for all $\xd \in (0, \xd_0]$ there holds
\be\la{4.25}
 \int_{K_{\xd}}
d^a |\nabla v| dx ~
\geq
~
 C \| d^b v \|_{L^{q}(K_{ \xd})},
 ~~~~~~\forall v \in
C_0^{\infty}(K_{\xd}).
\ee
\end{lemma}
The proof is quite similar to that of the previous Lemma.
The only difference is that instead of (\ref{12.5})
one uses (\ref{12.13}).  We omit the details.

We next have

\begin{theorem}
\la{thm4.1}
Let $\xO \subset \R^n$ be a  bounded  domain
 and $K$ a $C^2$ surface of co-dimension $k$, with $1 \leq k<n$,
satisfying condition ($R$).  We also assume
\be\la{condd}
 1 \leq p<n, \quad    p < q \leq \frac{pn}{n-p}, \quad {\rm and}
  \quad b=a-1 + \frac{q-p}{qp}n,
\ee
and set  $a = \frac{p-k}{p}$.
 Then
 there exists a   $\xd_0=\xd_0(p,q,\xO,K)$  and  $C= C(p,q,n,k)>0$ such that
for all $\xd \in (0, \xd_0]$
and  all $v \in C_0^{\infty}(\xO \setminus K)$ there holds
\be\la{4.30a}
 \int_{K_{\xd}}d^{p-k} |\nabla v|^p dx
+ \int_{\partial K_{\xd}} d^{1-k} |v|^p dS_x~
\geq~
C \|d^b v\|_{L^q(K_{\xd})}^p;
\ee
in particular  the constant  $C$ is independent of $\xO$, $K$.
\end{theorem}

{\em Proof:} We will use Lemma \ref{lem.41}.  Since in this Lemma the parameters
 $a$, $b$, $p$, $q$ have a different meaning,
to avoid confusion,
 we will use capital letters for the parameters
 $a$, $b$, $p$, $q$ appearing in the statement of the present  Theorem. That is,
we suppose that
\be\la{condb}
 1 \leq P<n, \quad    P < Q \leq \frac{Pn}{n-P}, \quad {\rm and}
  \quad B=A-1 + \frac{Q-P}{QP}n,
\ee
and for $A=\frac{P-k}{P}$, we will prove that the following estimate holds true
\be\la{4.30}
 \int_{K_{\xd}}d^{P-k} |\nabla v|^P dx
+ \int_{\partial K_{\xd}} d^{1-k} |v|^P dS_x~
\geq~
C \|d^B v\|_{L^Q(K_{\xd})}^P.
\ee

 We will argue in a similar
way,  as in the proof
of Theorem \ref{thm1}.
We first prove the following  $L^Q-L^P$ estimate:
\bea\la{4.31}
C(P,Q,n,k) \|d^B v\|^P_{L^Q(K_{\xd})}   &  \leq &
\int_{K_{\xd}}d^{P-k} |\nabla v|^P dx +
\int_{\partial K_{\xd}} d^{1-k} |v|^P dS_x  \nonumber \\
& & +
\|d^{\frac{P-k}{P}} v \|_{L^{\frac{(n-1)P}{n-P}}(\partial K_{\xd})}^P.
\eea

To this end   we replace in (\ref{4.17})  $v$ by $|v|^{s}$
with
\be\la{4.32}
s= Q\frac{P-1}{P} +1.
\ee
 Also, for $A$, $B$, $P$ and $Q$ as in (\ref{condb}),
 we set
\be\la{4.32a}
q= Qs^{-1},~~~~
b=Bs, ~~~~
 a= b+1 - \frac{q-1}{q}N = BQ\frac{P-1}{P}+A.
\ee
It is easy to check that $a$, $b$, $q$
thus defined satisfy (\ref{4.15}).
Then, from (\ref{4.17}) we have
\be\la{4.33}
 \|d^B v \|_{L^Q(K_{\xd})}^{1+\frac{P-1}{P}Q} = \| d^b |v|^s \|_{L^{q}(K_{\xd})}
 \leq    C \, s \,
\int_{K_{\xd}} d^a |v|^{s-1} |\nabla v| dx +
 C \int_{\partial K_{ \delta}} d^{a} |v|^s dx,
\ee
with $C= C(a,q,n,k)=C(P,Q,A,n,k)$.
Using Holder's inequality in the middle term of the right hand side we get
\bea\la{4.34}
\int_{K_{\xd}} d^a |v|^{s-1} |\nabla v| dx &  = &  \int_{K_{\xd}} d^{A} |\nabla v| ~~~
 d^{BQ \frac{P-1}{P}}
 |v|^{Q\frac{P-1}{P}} dx  \nonumber \\
& \leq &  \|d^A |\nabla v| \|_{L^P(K_{\xd})} ~~~ \|d^B v \|_{L^Q(K_{\xd})}^{\frac{P-1}{P}Q}, \nonumber \\
& \leq & c_{\xe}  \|d^A |\nabla v|
\|_{L^P(K_{\xd})}^{1+\frac{P-1}{P}Q} + \xe \|d^B v
\|_{L^Q(K_{\xd})}^{1+\frac{P-1}{P}Q}. \eea From now on we use the
specific value of $A=\frac{P-k}{P}$. For this choice of $A$
 a straightforward calculation shows that
\be\la{4.34b}
a-1+k = \frac{P-1}{P} \, \frac{Q-P}{P} \, (n-k) \neq 0,
\ee
and therefore it  corresponds to an acceptable value of  $a$,
see  (\ref{4.15}).   Because of (\ref{4.34b}) the case $k=n$ is excluded.

We next estimate the last term of (\ref{4.33}). Using Holder's
inequality (similarly as in Lemma \ref{interlem}), we get
\bea
 \int_{\partial K_{ \xd}} d^{a} |v|^s dx  \; = \; \int_{\partial K_{ \delta}} d^{\mu} |v|^{\xl (Q\frac{P-1}{P}+1)} ~~  d^{BQ \frac{P-1}{P}+A-\mu}
|v|^{(1-\xl)(Q\frac{P-1}{P}+1)} dx
  \nonumber \\
 \leq  \left( \int_{\partial K_{ \delta}}
d^{\frac{(P-k)(n-1)}{n-P}} |v|^{\frac{P(n-1)}{n-P}}dx  \right)^{ \frac{\xl(n-P)}{(n-1)P}(Q\frac{P-1}{P}+1)}
 \left( \int_{\partial K_{ \delta}} d^{1-k} |v|^P dx \right)^{\frac{1-\xl}{P} (Q\frac{P-1}{P}+1)},
 \nonumber
\eea
where,
\[
\xl=\frac{(n-1)(Q-P)}{Q(P-1)+P}, ~~~~~{\rm and}~~~~~ \mu=\frac{(n-1)(Q-P)(P-k)}{P^2}.
%~~~ r = \frac{P^2}{(n-P)(Q-P)}
\]
 Using then Young's inequality (cf (\ref{4.14})) we obtain for a positive
constant $C=C(P,Q,n)$,
\be\la{4.35}
 C\int_{\partial K_{ \xd}} d^{a} |v|^s dx
\leq  \left( \|d^{\frac{P-k}{P}} v
\|_{L^{\frac{P(n-1)}{n-P}}(\partial K_{\xd})} +
\|d^{\frac{1-k}{P}} v \|_{L^{P}(\partial K_{\xd})}
\right)^{Q\frac{P-1}{P}+1}. \ee From (\ref{4.33}), (\ref{4.34})
and (\ref{4.35}) we easily obtain (\ref{4.31}).

To complete the proof of the Theorem we  will show that
\be\la{4.37}
C \|d^{\frac{P-k}{P}} v \|^P_{L^{\frac{P(n-1)}{n-P}}(\partial K_{\xd})} \leq
\int_{K_{\xd}}d^{P-k} |\nabla v|^P dx +
\int_{\partial K_{\xd}} d^{1-k} |v|^P dS_x,
\ee
for a positive constant $C=C(P,Q,n,k)$. The proof of (\ref{4.37}) parallels that
of (\ref{12.19}). In particular, for $k=1$ this is precisely estimate (\ref{12.19}).
In the sequel we will sketch the proof of (\ref{4.37}).

Applying the critical trace inequality (\ref{12.21}) to $d^{\frac{P-k}{P}+\theta } v$,
$\theta>0$,
in the domain $K_{\xd}$ we obtain for $\xd$ sufficiently small the analogue of (\ref{12.23}),
that is
\bea\la{4.39}
\| d^{\frac{P-k}{P}} v\|^P_{L^{\frac{P(n-1)}{n-P}}(\partial K_{\xd})} &  \leq &
C(P,k)  \xd^{-\theta P} \int_{K_{\xd}} d^{P-k + P \theta} |\nabla v|^P dx  \nonumber \\
&  &  + \,\,
C(n,P,k,\theta) \, \,
  \xd^{-\theta P}  \int_{K_{\xd}} d^{-k + P \theta} |v|^P dx.
\eea We next estimate the last term of (\ref{4.39}). Starting from
the identity \be\la{4.40} (1-k +\theta P)  d^{-k + \theta P} =  -
d^{1-k + \theta P} \Delta d + {\rm div} (d^{1-k + \theta P} \nabla
d), \ee we multiply it by $|v|^P$ and  integrate by parts  over
$K_{\xd}$ to get \bea (1-k+\theta P ) \int_{K_{\xd}}  d^{-k +
\theta P} |v|^P dx = -\int_{K_{\xd}}
d^{1-k+\theta P} \Delta d |v|^P dx  \nonumber \\
   -P \int_{K_{\xd}} d^{1-k+\theta P} |v|^{P-1} \nabla d
 \cdot  \nabla |v|
 \, dx + \int_{\partial K_{\xd}} d^{1-k+\theta P} |v|^{P} dS_x,  \nonumber
\eea
or, equivalently,
\bea
\theta P  \int_{K_{\xd}}  d^{-k + \theta P} |v|^P dx = -\int_{K_{\xd}}
d^{k+\theta P} (d\Delta d+1-k) |v|^P dx  \nonumber \\
   -P \int_{K_{\xd}} d^{1-k+\theta P} |v|^{P-1} \nabla d
 \cdot  \nabla |v|
 \, dx + \int_{\partial K_{\xd}} d^{1-k+\theta P} |v|^{P} dS_x.  \nonumber
\eea
By our condition  ($R$) we have that $|d \Delta d +1-k| \leq c_0 d $.
On the other hand
\bea
 |P \int_{K_{\xd}} d^{1-k+\theta P} |v|^{P-1} \nabla d  \cdot  \nabla |v|
 dx | &   \leq  &  P \int_{K_{\xd}} d^{1-k+\theta P} |v|^{P-1} |\nabla v|  dx
\nonumber \\
 &   \leq  & P \e   \int_{K_{\xd}}  d^{-k + \theta P} |v|^P dx
+ P  c_{\e}  \int_{K_{\xd}} d^{P-k +  \theta P} |\nabla v|^P dx.
\nonumber
\eea
Putting together the last estimates we obtain, for $\e$, $\xd$ small
 the analogue of (\ref{12.25})
that is
\be\la{4.41}
C(P,\theta) \int_{K_{\xd}}  d^{-k +P \theta } |v|^P dx \leq
C(P)  \int_{K_{\xd}} d^{P-k + P \theta} |\nabla v|^P dx
 + \int_{\partial K_{\xd}} d^{1-k+P \theta } |v|^{P} dS_x.
\ee
Combining (\ref{4.39}), (\ref{4.41}) and using the fact that $d(x) \leq \xd$ when
$x \in K_{\xd}$, we complete the proof of (\ref{4.37}) as well as of the Theorem.

$\hfill \Box$

{\bf  Remark 1} We note that estimate (\ref{4.30a}) fails when $k=n$
 (see  (\ref{4.34b})).  This is not accidental as we shall see in the next
section.

{\bf  Remark 2} The  choice $a=\frac{p-k}{p}$  corresponds to
 the Hardy--Sobolev
inequality as it will become clear in the next section.
We note that the corresponding estimate
for  $a \in \R$
and $b$, $p$, $q$ as in (\ref{condd}) remains  true.
Thus, there exists a positive constant $C=C(a,n,p,q,k)$ such that
for all $v \in C^{\infty}_0(\xO \setminus K)$ there holds
\be\la{4.43}
\int_{K_{\xd}}  d^{ap}
 |\nabla v |^p dx +
\int_{  \partial K_{\xd}}
 d^{(a-1)p+1}  | v|^p dS_x
\geq C \, \|d^b v\|_{L^q(\xO)}.
\ee
The proof of (\ref{4.43}) in case $a \neq \frac{p-k}{p}$ is much simpler than
in the case $a=\frac{p-k}{p}$.
We also note that if  $a \neq \frac{p-k}{p}$ then (\ref{4.43}) is true even if  $k=n$.

We will finally prove the analogue of Theorem \ref{thm3}.

\begin{theorem}
\la{thm4.2}
Let $\xO \subset \R^n$ be a domain  and $K$ a  surface of codimension $k$,
 $1 \leq k<n$,
satisfying both conditions ($R$). In addition we assume that
 $D = \sup_{x \in \xO}d(x)< \infty$,
condition $\C$ is satisfied  and \be\la{cond2}
 1 \leq p<n, \quad    p < q \leq \frac{pn}{n-p}, \quad {\rm and}
  \quad b=a-1 + \frac{q-p}{qp}n.
\ee
We set  $a = \frac{p-k}{p}$.
 Then
there exists a positive constant $C=C(p,n,\xO,K)$ such that
for all $v \in C_0^{\infty}(\xO  \setminus K)$ there holds
\be\la{4.50a}
\int_{\xO}d^{p-k} |\nabla v|^p dx +
|\int_{\xO} d^{-k}(-d \Delta d -1+k) |v|^p dx|
\geq
C \|d^b v\|_{L^q(\xO)}^p.
\ee
\end{theorem}
{\em Proof:} As before,
to avoid confusion in the proof, we will use capital letters for the parameters
 $a$, $b$, $p$, $q$ appearing in the statement of the Theorem. That is,
we suppose that
\[
 1 \leq P<n, \quad    P < Q \leq \frac{Pn}{n-P}, \quad {\rm and}
  \quad B=A-1 + \frac{Q-P}{QP}n,
\]
and for $A=\frac{P-k}{P}$, we will prove that
\be\la{4.50}
\int_{\xO}d^{P-k} |\nabla v|^P dx +
|\int_{\xO} d^{-k}(-d \Delta d -1+k) |v|^P dx|
\geq
C \|d^B v\|_{L^Q(\xO)}^P.
\ee

Let $\xa(t) \in C^{\infty}([0,\infty))$ be   the  nondecreasing function
defined at the beginning of the proof of Theorem \ref{thm2} and
 $\phi_{\xd}(x):=\xa(\frac{d(x)}{\xd}) \in C^2_0(\xO)$, so that
 $\phi_{\xd}=1$ on $K_{\delta/2}$,
 $\phi_{\xd}=0$ on $K_{ \delta}^c$ and $| \nabla \phi_{\xd}| \leq  \frac{C_0}{\xd}$
with  $C_0$ a universal constant.

For $v \in  C_0^{\infty}(\xO)$ we write  $v= \phi_{\xd} v + (1-\phi_{\xd}) v$.
The function $\phi_{\xd} v$ is compactly supported in $K_{\xd}$, and by Lemma \ref{lem.42} we
have
\be\la{4.51}
 C(a,n,q) \| d^b v \|_{L^{q}(K_{ \xd})} \leq   \int_{K_{\xd}}
d^a |\nabla v| dx.
\ee
On the other hand $(1-\phi_{\xd}) v$  is compactly supported in $K_{\xd/2}^c$ and using
(\ref{12.13}) we easily get
\be\la{4.53}
\|d^b (1-\phi_{\xd}) v\|_{L^q(K_{\xd/2}^c)} \leq C(\xO)
 \frac{D^{|b|}}{\xd^{|a|}}
 \|d^a |\nabla ((1-\phi_{\xd}) v)| \|_{L^1(K_{ \xd/2}^c)}.
\ee
Combining (\ref{4.51}) and (\ref{4.53}) we obtain   the analogue of (\ref{12.45})
which is
\be\la{4.55}
C \| d^a v \|_{L^{\frac{n}{n-1}}(\xO)} \leq
 \int_{\xO}  | d^a  \nabla v | dx +
\int_{K_{ \delta} \setminus K_{ \delta/2}} d^{a-1}  |v| dx.
\ee
We next pass to  $L^Q$--$L^P$ estimates.  We replace in (\ref{4.55})  $v$ by $|v|^{s}$
with $s$ as in (\ref{4.32}). Also, for $A=\frac{P-k}{P}$
and $B$, $P$, $Q$ as in
(\ref{4.32a}), we get (cf (\ref{4.33}))
\be\la{4.57}
 C \|d^B v \|_{L^Q(K_{\xd})}^{1+\frac{P-1}{P}Q}
%\| d^b |v|^s \|_{L^{q}(K_{\xd})}
 \leq    s \,
\int_{K_{\xd}} d^a |v|^{s-1} |\nabla v| dx +
  \int_{K_{ \delta} \setminus K_{ \delta/2} } d^{a-1} |v|^s dx.
\ee
Using Holder's inequality in both terms of the right hand side we get
\bea
\int_{\xO} d^a |v|^{s-1} |\nabla v| dx &  = &  \int_{\xO} d^{A} |\nabla v| ~~~
 d^{BQ \frac{P-1}{P}}
 |v|^{Q\frac{P-1}{P}} dx  \nonumber \\
& \leq &  \|d^A |\nabla v| \|_{L^P(\xO)} ~~~ \|d^B v \|_{L^Q(\xO)}^{\frac{P-1}{P}Q},
\nonumber
\eea
and
\bea
 \int_{K_{ \delta} \setminus K_{ \delta/2}} d^{a-1} |v|^s dx & = &
 \int_{K_{ \delta} \setminus K_{ \delta/2}} d^{A-1} |v| ~~~
 d^{BQ \frac{P-1}{P}}
 |v|^{Q\frac{P-1}{P}} dx  \nonumber \\
& \leq &   \|d^{A-1} | v| \|_{L^P(K_{ \delta} \setminus K_{
\delta/2})} ~~~ \|d^B v \|_{L^Q(\xO)}^{\frac{P-1}{P}Q}. \nonumber
\eea Substituting into (\ref{4.57}) we get after simplifying,
\be\la{4.59}
 C \|d^B v \|^P_{L^Q(\xO)} \leq
 \int_{\xO} d^{P-k}  |\nabla v|^P dx  +
  \int_{K_{ \delta} \setminus K_{ \delta/2}}
  d^{-k}  |v|^P dx.
\ee
Here we have also used the specific value of $A=\frac{P-k}{P}$.
To conclude we need to estimate the last term in (\ref{4.59}).
For $\theta>0$, we clearly have
\be\la{4.61}
\left( \frac{\xd}{2} \right)^{p \theta}
 \int_{K_{ \delta} \setminus K_{ \delta/2}}  d^{-k} |v|^P dx
\leq
\int_{K_{ \delta} \setminus K_{ \delta/2}}  d^{-k+P \theta} |v|^P dx
\leq
\int_{\xO} d^{-k+P \theta} |v|^P dx.
\ee
To estimate the last term we work as in
%(\ref{12.23a})--(\ref{12.25})
 (\ref{12.50})--(\ref{12.51})
(see also (\ref{4.40})--(\ref{4.41}))
 to finally get
\be\la{4.63}
\int_{\xO} d^{-k+P \theta} |v|^P dx \leq
C(p)  \int_{\xO}      d^{P-k + P \theta} |\nabla v|^P dx
 + | \int_{\xO}   d^{-k+ P \theta }(-d \Delta d +1-k) |v|^{P} dx| .
\ee We not that we also used the fact that \be\la{4.64} p\neq k,
~~~{\rm and }~~~~(p-k)(d \Delta d +1-k) \leq 0, ~~{\rm on} ~~
\Omega \setminus K, \ee which is a direct consequence of condition
$\C$;  see \cite{BFT}. Combining (\ref{4.61}) and (\ref{4.63}) and
recalling that $d \leq D$, we get \be\la{4.65}
 C(P, \theta, \frac{\xd}{D})\,
 \int_{K_{ \delta} \setminus K_{ \delta/2}}  d^{-k} |v|^P dx
\leq
 \int_{\xO}       d^{P-k} |\nabla v|^P dx
+| \int_{\xO}   d^{-k }(-d \Delta d +1-k)) |v|^{P} dx|,
\ee
and  the result  follows easily.

$\hfill \Box$
\\
{\bf Remark 1} As in Theorem \ref{thm4.1} the case $k=n$ is excluded.
\\
{\bf Remark 2} In case $k=1$ or in case $q =\frac{np}{n-p}$,
the dependence of  the constant
$C$ in (\ref{4.50a}) from $\xO$, $K$ is the same as in Theorem \ref{thm3},
that is, $C=C(n,p,q,c_0D)$.
\\
{\bf Remark 3 } In case $a \neq \frac{p-k}{p}$ the analogue of (\ref{4.50a}) remains true.
That is, for $b$, $p$, $q$ as in (\ref{cond2})
\be\la{4.80}
\int_{\xO}d^{ap} |\nabla v|^p dx +
|\int_{\xO} d^{(a-1)p}(-d \Delta d -1+k) |v|^p dx|
\geq
 C \, \|d^b v\|_{L^q(\xO)},
\ee
for a constant $C=C(p,q,n,k,a)>0$. The case $k=n$ is not excluded.

\setcounter{equation}{0}
\section{Extended Hardy--Sobolev inequalities}
In this Section we will use the $v$--inequalities of the previous Section to
prove new Hardy--Sobolev inequalities.
For $V \subset \R^n$ we set
\be\la{5.1}
I_{p,k}[u](V)  :\int_{V} |\nabla u|^p dx- \left| \frac{p-k}{p}\right|^{p} \int_{V} \frac{|u|^p}{d^p} dx.
\ee
Then for $u(x)=d^{H}(x) v(x)$ with
\[
H:=\frac{p-k}{p},
\]
  we have for $p \geq 2$,
\be\la{5.3}
I_{p,k}[u](V) \geq c(p) \int_{V}  d^{p-k} |\nabla v|^p dx +
H |H|^{p-2}  \int_{ V} d^{1-k}  \nabla d \cdot \nabla |v|^pdx.
\ee
The proof of (\ref{5.3}) is quite  similar to the proof of  (\ref{3.5}).

As in the previous section,
\be\la{cond5}
 1 \leq p<n, \quad    p < q \leq \frac{pn}{n-p}, \quad {\rm and}
  \quad b=a-1 + \frac{q-p}{qp}n.
\ee
We will be interested in the specific value  $a= \frac{p-k}{p}$ which corresponds to
the critical Hardy Sobolev inequalities.

We first present estimates in $K_{\xd}$.
\begin{theorem}\la{thm5.1}
Let $2 \leq p <n$  and  $p < q \leq \frac{np}{n-p}$.
We assume that  $\xO \subset \R^n$ is  a  bounded  domain
 and $K$ a $C^2$ surface of co-dimension $k$, with $1 \leq k<n$,
satisfying condition ($R$).
 Then, there exist  positive constants $C=C(n,k,p,q)$  and  $\xd_0=\xd_0(p,n,\xO,K)$
such that  for $0< \xd \leq \xd_0 $  and   $ u \in C_0^{\infty}(\xO \setminus K)$   we have: \\
(a) If  $p>k$ then
\be\la{5.10}
\int_{K_{\xd}} |\nabla u|^p dx- |H|^p
 \int_{K_{\xd}} \frac{|u|^p}{d^p} dx
\geq C
 \left( \int_{K_{\xd}}  d^{-q+\frac{q-p}{p}n} |u|^{q} dx \right)^{\frac{p}{q}}.
\ee
(b) If $p <k$, the Hardy inequality
\be\la{5.10b}
\int_{K_{\xd}} |\nabla u|^p dx- |H|^p
 \int_{K_{\xd}} \frac{|u|^p}{d^p} dx
\geq 0,
\ee
in general fails.
 However, there exists  a positive
constant $M$ such that
\be\la{5.10c}
\int_{K_{\xd}} |\nabla u|^p dx- |H|^p
 \int_{K_{\xd}} \frac{|u|^p}{d^p} dx + M \int_{K_{\xd}} | u|^p dx
\geq C
 \left( \int_{K_{\xd}}  d^{-q+\frac{q-p}{p}n} |u|^{q} dx \right)^{\frac{p}{q}}.
\ee
We emphasize  that   $C=C(n,k,p,q)>0$ is independent of $\xO$, $K$.  \\
 (c)  If in
addition,  $u$ is supported in $K_{\xd}$, that is
$u \in  C_0^{\infty}(K_{\xd} \setminus K)$
then, (\ref{5.10}) holds true even for $p<k$.
\end{theorem}

{\em Proof:}  Using    (\ref{5.1}) and  integrating  by parts  once   we have that
\bea\la{5.11}
 I_{p,k}[u](K_{\xd}) &  \geq &  C(p)  \int_{K_{\xd}}  d^{p-k} |\nabla v|^p  dx +
 H  |H|^{p-2} \int_{K_{\xd}}   d^{-k }(-d \Delta d +k-1)) |v|^{p} dx  \nonumber \\
 &   &  +  ~ H  |H|^{p-2} \int_{\partial K_{\xd}} d^{1-k} |v|^p dS_x.
\eea
At first we  estimate the middle term of the right hand side.
We have that
\be\la{5.12}
|d \Delta d +1-k| \leq c_0 d,~~~~~{\rm   for} ~~~~ x \in K_{\xd},
\ee
 and therefore
\be\la{5.13}
| \int_{K_{\xd}}  d^{-k }(-d \Delta d +k-1)) |v|^p dx| \leq c_0 \int_{K_{\xd}}d^{1-k}  |v|^p dx.
\ee

At this point we will derive some general estimates that we  will use in the sequel.
Our goal is to prove (\ref{5.15a}) and (\ref{5.15b})  below.
For $a \in \R$ we consider  the identity
 $ (1+a)  d^{a} + d^{1+a} \Delta d = {\rm div}(d^{1+a} \nabla d)$. Multiply
by $|v|^p$ and integrate
 by parts to get
\[
 (a+1) \int_{K_{\xd}}  d^{a} |v|^p dx +  \int_{K_{\xd}} d^{a+1}
 \Delta d |v|^P dx  \; = \; - p \int_{K_{\xd}} d^{a+1} \nabla d
 \cdot \nabla |v| \,  |v|^{p-1}  dx + \int_{\partial K_{\xd}} d^{a+1}
|v|^p dS_x,
\]
or, equivalently,
\bea\la{5.14}
 (a+k) \int_{K_{\xd}}  d^{a} |v|^p dx +  \int_{K_{\xd}} d^{a}
 (d \Delta d + 1-k) |v|^p dx  =        \nonumber \\
 - p \int_{K_{\xd}} d^{a+1} \nabla d
 \cdot \nabla |v| \,  |v|^{p-1}  dx + \int_{\partial K_{\xd}} d^{a+1}
|v|^p dS_x.
\eea
We next estimate the first term of
the right hand side of (\ref{5.14})
\bea\la{5.14a}
p\int_{K_{\xd}} d^{a+1} \nabla d
 \cdot \nabla |v| \,  |v|^{p-1}  dx &  \leq  &
\left( \int_{K_{\xd}} d^{a} | v|^p dx \right)^{\frac{p-1}{p}}
\left( \int_{K_{\xd}} d^{a+p} |\nabla v|^p dx \right)^{\frac1p}
  \nonumber \\
&  \leq  &  \xe (p-1)  \int_{K_{\xd}} d^{a} | v|^p dx +
 \xe^{-(p-1)} \int_{K_{\xd}}
 d^{a+p} |\nabla v|^p dx.
\nonumber \eea From  this,  (\ref{5.12}) and (\ref{5.14}) we
easily obtain the following two estimates: \be\la{5.15a} (|a+k|
-c_0 \xd -\xe(p-1)) \int_{K_{\xd}}  d^{a} |v|^p dx \leq
 \xe^{-(p-1)} \int_{K_{\xd}}
 d^{a+p} |\nabla v|^p dx + \int_{\partial K_{\xd}} d^{a+1}
|v|^p dS_x,
\ee
and,
\be\la{5.15b}
 \int_{\partial K_{\xd}} d^{a+1}
|v|^p dS_x  \leq \xe^{-(p-1)} \int_{K_{\xd}}
 d^{a+p} |\nabla v|^p dx +
(|a+k| +c_0 \xd +\xe(p-1)) \int_{K_{\xd}}  d^{a} |v|^p dx. \ee
From (\ref{5.14a}) taking $a=1-k$ we get that \be\la{5.16}
\int_{K_{ \xd}}d^{1-k} |v|^p dx  \leq C(p)  \xd \int_{K_{
\xd}}d^{p-k} |\nabla v|^p dx +  C(p) \xd \int_{\partial
K_{\xd}}d^{1-k} |u|^p  dS_x. \ee

At this point we distinguish two cases according to whether $p>k$ or $p<k$.
Assume first that  $p>k$, or equivalently, $H >0$. Then from  (\ref{5.11}) and  (\ref{5.16})
we get that
\be\la{5.17a}
 I_{p,k}[u](K_{\xd})   \geq  C(p)  \int_{K_{\xd}}  d^{p-k} |\nabla v|^p  dx
 + C(p,k)  \int_{\partial K_{\xd}} d^{1-k} |v|^p dS_x.
\ee
Using  Theorem \ref{thm4.1} as well as the fact that
\[
\|d^b v\|_{L^q(K_{\xd})}^p \; = \; \left( \int_{K_{\xd}}
d^{-q+\frac{q-p}{p}n} |u|^{q} dx \right)^{\frac{p}{q}},
\]
we easily  obtain (\ref{5.10}).

If $u \in  C_0^{\infty}(K_{\xd}\setminus K)$
then the boundary terms in (\ref{5.11}) and (\ref{5.16})
 are  absent and the same argument
yields (\ref{5.10}) even if  $p<k$.

Suppose now that  $p<k$, that is,   $H<0$.
Using   again  (\ref{5.11}) and  (\ref{5.16}) we get that
\be\la{5.17}
 I_{p,k}[u](K_{\xd})   \geq   C(p)  \int_{K_{\xd}}  d^{p-k} |\nabla v|^p  dx
 - C(p,k)  \int_{\partial K_{\xd}} d^{1-k} |v|^p dS_x.
\ee
To estimate the last term of this we will use (\ref{5.15b}) with $a=p-k$
in the following way
\bea\la{5.19}
\int_{\partial K_{\xd}} d^{1-k} |v|^p dS_x  &  =   & \xd^{-p}
\int_{\partial K_{\xd}} d^{1+p-k} |v|^p dS_x    \nonumber \\
& \leq &
\xe^{-(p-1)} \int_{K_{\xd}}
 d^{p-k} |\nabla v|^p dx + C(\xe,p)  \xd^{-p}
\int_{K_{\xd}}  d^{p-k} |v|^p dx. \eea From  (\ref{5.17}) and
(\ref{5.19}) choosing $\xe$ big we get \be\la{5.19a}
 I_{p,k}[u](K_{\xd})   \geq   C(p)  \int_{K_{\xd}}  d^{p-k} |\nabla v|^p  dx
 - M  \int_{ K_{\xd}} d^{p-k} |v|^p dx.
\ee
On the other hand from (\ref{5.19}) and Theorem 4.4 we get that
\be\la{5.19b}
C(p,q,n,k) \|d^b v\|_{L^q(K_{\xd})}^p  \leq
 C(p) \int_{K_{\xd}}  d^{p-k} |\nabla v|^p  dx +
M \int_{ K_{\xd}} d^{p-k} |v|^p dx. \ee From (\ref{5.19a})  and
(\ref{5.19b}) we easily conclude (\ref{5.10c}).

It remains to explain why    when $p <k$ and  $ u \in C_0^{\infty}(\xO \setminus K)$
 the simple Hardy  (\ref{5.10b}) in general  fails. Let us consider the case where  $K$ and therefore
$K_{\xd}$ are strictly contained in $\xO$. In this case the function
 $u_{\xe}=d^{H+\xe}$, for $\xe>0$  is in $W^{1,p}(K_{\xd})$. On the other hand for $p <k$
a simple density argument shows that $W^{1,p}(K_{\xd} \setminus K)=W^{1,p}(K_{\xd})$.
 An easy
calculation shows that
\be\la{5.20}
\int_{K_{\xd}} |\nabla u_{\xe}|^p dx- |H|^p
 \int_{K_{\xd}} \frac{|u_{\xe}|^p}{d^p} dx = (|H+\xe|^p-|H|^p) \int_{K_{\xd}} d^{-k+p\xe}dx<0,
\ee
by taking $\xe>0$ small and noting that $H<0$.

$\hfill \Box$
\\
{\bf Remark}  The result is not true in case $k=n$, as discussed in the
introduction.

We next  prove  estimates in $\xO$.
\begin{theorem}\la{thm5.2}
Let $2 \leq p <n$  and  $p < q \leq \frac{np}{n-p}$.
We assume that  $\xO \subset \R^n$ is  a  bounded  domain
 and $K$ a $C^2$ surface of co-dimension $k$, with $1 \leq k<n$,
satisfying condition ($R$).
Then,
 there exist  positive constants $C=C(n,k,p,q)$ and  $M$   such that
for all $ u \in C_0^{\infty}(\xO \setminus K)$,
there holds
\be\la{5.21}
\int_{\xO} |\nabla u|^p dx-  \left| \frac{p-k}{p} \right|^p  \int_{\xO} \frac{|u|^p}{d^p} dx
 + M   \int_{\xO} |u|^p dx   \geq
 C \left( \int_{\xO}  d^{-q+\frac{q-p}{p}n} |u|^{q} dx \right)^{\frac{p}{q}},
\ee
We note that $C(n,k,p,q)$ is independent of $\xO$, $K$.
\end{theorem}

{\em Proof:} Clearly we have
\be\la{5.23}
I_{p,k}[u](\xO)=I_{p,k}[u](K_{\xd})+ I_{p,k}[u](K_{\xd}^c).
\ee
By Theorem    \ref{thm5.1}    for $\xd$ small we have
\be\la{5.25}
I_{p,k}[u](K_{\xd}) \geq C(n,k,p,q)
 \left( \int_{K_{\xd}}  d^{-q+\frac{q-p}{p}n} |u|^{q} dx \right)^{\frac{p}{q}}
- M \int_{K_{\xd}}|u|^{p} dx.
\ee
Since $d(x) \geq \xd$ in $K_{\xd}^c$,
\be\la{5.27}
I_{p,k}[u](K_{\xd}^c) \geq \int_{K_{\xd}^c} |\nabla u|^p dx -
 %\left( \frac{p-1}{p \xd} \right)^p
C(p,k,\xd)
 \int_{K_{\xd}^c} |u|^p dx.
\ee From  the Sobolev embedding of $L^{\frac{np}{n-p}}(K_{\xd}^c)$
into
 $W^{1,p}(K_{\xd}^c)$ we get
\[
 \| u \|^p_{L^{\frac{np}{n-p}}(K_{\xd}^c)}  \leq
C(p,n) \int_{K_{\xd}^c} |\nabla u|^p dx + C(p,n, \xO,K) \int_{K_{\xd}^c} |u|^p dx.
\]
Using  the interpolation Lemma \ref{interlem} (with $a=0$) we have
\bea\la{5.29}
 C(n,p,q) \, \left( \int_{K_{\xd}^c}
 d^{-q+\frac{q-p}{p}n} |u|^{q} dx \right)^{\frac{p}{q}}
& \leq &
 \| u\|^p_{L^{\frac{pn}{n-p}}(K_{\xd}^c)} +
\|d^{-1}u\|^p_{L^p(K_{\xd}^c)},  \nonumber \\
& \leq &
 \| u\|^p_{L^{\frac{pn}{n-p}}(K_{\xd}^c)} +
  \xd^{-p} \|u\|^p_{L^p(K_{\xd}^c)}.
\eea From  (\ref{5.27})--(\ref{5.29}) we get for $M=M(n,p,q,
\xO,K)$, \be\la{5.31} I_{p,k}[u](K_{\xd}^c) \geq C(n,p,q) \left(
\int_{K_{\xd}^c}
 d^{-q+\frac{q-p}{p}n} |u|^{q} dx \right)^{\frac{p}{q}}
- M  \int_{K_{\xd}^c} |u|^p dx.
\ee
The result follows from (\ref{5.23}), (\ref{5.25}) and (\ref{5.31}).

$\hfill \Box$

Our final result reads:
\begin{theorem}\la{thm5.3}
Let $2 \leq p <n$  and  $p < q \leq \frac{np}{n-p}$.
We assume that   $\xO \subset  \R^n$
is  a domain  and $K$ a  surface of co-dimension $k$,
 $1 \leq k<n$,
satisfying condition ($R$). In addition we assume that
 $D = \sup_{x \in \xO}d(x)< \infty$ and
condition $\C$ is satisfied.
 Then
 for all $ u \in C_0^{\infty}(\xO)$
there holds
\be\la{5.35}
\int_{\xO} |\nabla u|^p dx- \left| \frac{p-k}{p}\right|^{p} \int_{\xO} \frac{|u|^p}{d^p} dx
\geq
 C \left( \int_{\xO}  d^{-q+\frac{q-p}{p}n} |u|^{q} dx \right)^{\frac{p}{q}},
\ee
for $C=C(n,P,Q,\xO,K)>0$.
\end{theorem}

{\em Proof:} Working as in the derivation of (\ref{5.11}) we get
\bea\la{5.37}
C(p,k) \, I_{p,k}[u](\xO) &  \geq &    \int_{\xO}  d^{p-k} |\nabla v|^p  dx +
 H  \int_{\xO}   d^{-k }(-d \Delta d +1-k)) |v|^{p} dx.
\eea Because of condition $\C$ we have that $H (-d \Delta d +1-k)
\geq 0$, see (\ref{4.64}). The result  then follows  from Theorem
\ref{thm4.2}.

 $\hfill \Box$

\end{document}